\documentclass[12pt,a4paper]{amsart}
\usepackage{mathtools}
\usepackage{graphicx}
\usepackage{amsthm}
\usepackage{amsfonts}
\usepackage{amssymb}
\usepackage{url}
\usepackage{enumerate}
\usepackage{geometry}
\usepackage{subcaption}
\usepackage{fullpage}

\title{Rigorous computation of Maass cusp forms of squarefree level}
\author{Andrei Seymour-Howell}
\address{School of Mathematics, Fry Building, University of Bristol, Woodland Road, Bristol, BS8 1UG, UK}
\email{andrei.seymour-howell@bristol.ac.uk}

\date{}

\newtheorem{theorem}{Theorem}[section]

\theoremstyle{remark}
\newtheorem*{remark}{Remark}

\theoremstyle{definition}

\newcommand{\SL}{\textup{SL}}

\newcommand{\NN}{\mathbb{N}}
\newcommand{\ZZ}{\mathbb{Z}}
\newcommand{\QQ}{\mathbb{Q}}
\newcommand{\RR}{\mathbb{R}}
\newcommand{\CC}{\mathbb{C}}
\newcommand{\HH}{\mathbb{H}}
\newcommand{\Real}{\textup{Re}}
\newcommand{\Imag}{\textup{Im}}

\newcommand{\sign}{\textup{sign}}
\newcommand{\sinc}{\textup{sinc}}

\begin{document}
\begin{abstract}
    We derive an algorithm to rigorously compute and verify Maass cusp forms of squarefree level and trivial character. The main tool we use is an explicit version of the Selberg trace formula with Hecke operators due to Str\"{o}mbergsson. We use this algorithm to compute several thousand Maass forms for a range of levels and use this data to obtain numerical evidence towards various conjectures.
    \end{abstract}
\maketitle
\section{Introduction}
In the 1990s Hejhal \cite{Hejhal} derived an algorithm to compute the Laplace eigenvalues and Fourier coefficients of Maass cusp forms for certain subgroups of $\SL_2(\ZZ)$. In 2006 Str\"omberg, in his PhD thesis \cite{stromberg}, further generalised this method to compute the eigenvalues of Maass cusp forms for congruence subgroups $\Gamma_0(N)$ with a character. Whilst this algorithm is numerically stable in practice, it unfortunately relies on a heuristic argument and has not yet been proven to rigorously converge to true Maass cusp forms.

There has been progress towards numerically verifying numerical computations of Maass cusp forms, most notably from Booker, Str\"{o}mbergsson and Venkatesh \cite{bookerstromvenka}, who derived a method to numerically verify Maass cusp forms for $\SL_2(\ZZ)$. Using this method they verified the first $10$ Laplace eigenvalues to $100$ decimal places. For general level $N$, there currently does not exist an algorithm to rigorously verify numerical computations of Maass cusp forms for $\Gamma_0(N)$.

In this paper, we derive an algorithm to numerically compute and rigorously verify the Laplace and Hecke eigenvalues for Maass cusp forms for $\Gamma_0(N)$ with $N$ squarefree and trivial character. The main tool that we use is an explicit version of the Selberg trace formula with Hecke operators derived by Str\"{o}mbergsson in \cite{Andreaspre}. The inclusion of the Hecke operators allows us to use linear algebra to help pick out individual eigenvalues.

The Selberg trace formula has been used for numerical computation before by Booker and Str\"{o}mbergsson in \cite{stfz} to numerically verify the Selberg eigenvalue conjecture. However those authors were focused on proving the non-existence of Maass forms in an interval, rather than computing individual examples. In the case of holomorphic modular forms, explicit versions of the Selberg trace formula have been used to compute bases of cusp forms, for example in \cite{modformspari}.

The outline of the paper is as follows. In Section \ref{sec:algorithm} we present the algorithm for computing and then rigorously verifying the Laplace eigenvalues and Hecke eigenvalues of Maass cusp forms for squarefree level and trivial character. In Section \ref{sec:traceformula} we state the explicit form of the Selberg trace formula that we use and explain computational aspects on how to compute it. In Section \ref{sec:testfunc} we choose and optimise the test function for the trace formula such that it maximises the precision of the computation. Finally, in Section \ref{sec:compres} we state the computational results and show some numerical evidence towards the Ramanujan--Petersson conjecture, Sato--Tate conjecture and the Riemann hypothesis for L-functions of Maass cusp forms.

\textbf{Acknowledgments.} I would like to thank my PhD supervisor Andrew Booker for his guidance and comments with this research and Andrew Sutherland for their comments on this paper.
\section{Algorithm}
\label{sec:algorithm}
In this section, we derive the algorithm to compute and rigorously verify the Laplace and Hecke eigenvalues for Maass cusp forms of squarefree level $N$. The central tool used here is the Selberg trace formula with Hecke operators. The main idea here is to use linear algebra to remove the contribution of all the forms up to some limit and isolate just one form. We then use our approximation to this form to see how well it removes the remaining contribution.
\subsection{Preliminary}
Let $\HH = \{ z = x+iy \in \CC \mid y > 0\}$ denote the hyperbolic upper half-plane. This is acted on by certain subgroups of the modular group $\SL_2(\ZZ)$ via linear fractional transformations. For $N \in \NN$, the subgroups we shall consider in this paper are the Hecke congruence subgroups $\Gamma_0(N) \subseteq \SL_2(\ZZ)$ defined by
\begin{align*}
    \Gamma_0(N) = \left\{ \left. \begin{pmatrix}
    a&b\\
    c&d
    \end{pmatrix} \in \SL_2(\ZZ) \right| c \equiv 0 \quad \textup{ mod } N\right\}.
\end{align*}
A Maass cusp form for $\Gamma_0(N)$ of weight $0$ is a non-constant, smooth function $f: \Gamma_0(N) \backslash \HH \to \CC$ that satisfies the following properties:
\begin{enumerate}
    \item $f(\gamma z) = f(z)$ for all $z \in \HH$ and $\gamma \in \Gamma_0(N)$;
    \item $f$ vanishes at the cusps of $\Gamma_0(N)$;
    \item $f \in L^2(\Gamma_0(N) \backslash \HH)$;
    \item $f$ is an eigenfunction of the Laplace--Beltrami operator $\Delta$ on $\HH$ given by
    \begin{align*}
        \Delta = -y^2 \left( \frac{\partial^2}{\partial x^2} + \frac{\partial^2}{\partial y^2} \right).
    \end{align*}
\end{enumerate}

Now if $K \mid N$, then $\Gamma_0(N) \subseteq \Gamma_0(K)$. Notably, if $f$ is a Maass cusp form of $\Gamma_0(K)$, then $f(kz)$ is a Maass cusp form of $\Gamma_0(N)$ for all $k \mid \frac{N}{K}$. Forms that arise like this for $\Gamma_0(N)$ we call ``oldforms''. We define ``newforms'' to be forms in the orthogonal complement (with respect to the Petersson inner product) of the space spanned by the oldforms.

In addition, we define the Hecke eigenvalues $a(n)$ of a level $N$ Maass form $f$ by the relation
\begin{align*}
    a(n) f = T_n f,
\end{align*}
where $T_n$ is the $n$th Hecke operator, defined by
\begin{align*}
    T_n f(z) = \frac{1}{\sqrt{|n|}} \underset{\substack{ad=n\\ (a,N) = 1 \\ d > 0}}{\sum} \sum_{j = 0}^{d-1} 
        \begin{dcases}
        f \left( \frac{az + j}{d} \right) & \text{if } n > 0,\\
        f \left( \frac{a \overline{z} + j}{d} \right) & \text{if } n < 0.
        \end{dcases}
\end{align*}

For newforms, we normalise the Hecke eigenvalues by setting $a(1)=1$. Furthermore, we say a Maass form is even if $a(-n) = a(n)$ and odd if $a(-n) = -a(n)$. 
More information about the theory of Maass cusp forms can be found in Bump \cite{bump} and Iwaniec \cite{iwaniec}. \\

Consider the space of Maass newforms of level N, Laplace eigenvalue $\lambda$ and trivial character, denoted by $S_{\lambda}(N)$. Let $\lbrace  f_j \rbrace_{j=1}^{\infty}$ be a sequence of normalised Hecke eigenforms such that it is a basis for $\bigoplus_{\lambda > 0} S_{\lambda}(N)$. Let $\lambda_j$ denote the Laplace eigenvalue of $f_j$ and assume that $\lambda_1 \leq \lambda_2 \leq \ldots$. In addition, let $a_j(n)$ be the Hecke eigenvalues for $f_j$.

The Selberg trace formula allows us to compute
\begin{align*}
    t(n,H) := \sum_{j=1}^{\infty} a_j(n) H(\lambda_j),
\end{align*}
for any non-zero $n \in \ZZ$ with $(n,N) = 1$ and any sufficiently nice test function $H$. Using the Hecke relations, we compute that
\begin{align*}
    \left( \sum_{m=1}^{M} c(m) a_j(m) \right)^2 = \sum_{m_1 = 1}^{M} \sum_{m_2 = 1}^{M} c(m_1) c(m_2) \sum_{d | (m_1,m_2)} a_j \left( \frac{m_1 m_2}{d^2} \right),
\end{align*}
for any sequence $\lbrace c(m) \rbrace_{m=1}^{M}$ of real numbers, satisfying $c(m) = 0$ whenever $(m,N) > 1$. Thus, defining
\begin{align*}
    Q(c,H) := \sum_{j=1}^{\infty} \left( \sum_{m=1}^{M} c(m) a_j(m) \right)^2 H(\lambda_j),
\end{align*}
we have 
\begin{align}
     Q(c,H) = \sum_{m_1 = 1}^{M} \sum_{m_2 = 1}^{M} c(m_1) c(m_2) \sum_{d | (m_1,m_2)} t \left( \frac{m_1 m_2}{d^2}, H \right).
     \label{eqn:Qeqn}
\end{align}

\subsection{Computing the forms}
Let $H$ be a non-negative test function and let $\widetilde{H}(\lambda) = \lambda H(\lambda)$. Let $Q$ and $\widetilde{Q}$ denote the respective matrices of the quadratic forms $Q(c,H)$ and $Q(c,\widetilde{H})$. Thus the approximations of the Laplace eigenvalues are the solutions to the generalised symmetric eigenvalue equation
\begin{align}
	\widetilde{Q}x = \lambda Q x.
	\label{eqn:symeig}
\end{align}
We solve this by first diagonalising $Q = PDP^T$, where $P$ is an orthogonal matrix and $D$ is diagonal with positive entries. Then the solutions to \eqref{eqn:symeig} will just be the eigenvalues of $D^{-1/2} P^T \widetilde{Q} P D^{-1/2}$. For each eigenvalue $\widetilde{\lambda}_i$, we set $c_i$ to be the corresponding eigenvector. We will use the components of $c_i$ to form the sequence $c(m)$ for each eigenvalue. 

\subsection{Verifying the forms}
Firstly, for the verification we shall prove that there exists a Laplace eigenvalue near $\widetilde{\lambda}_i$. For this, we define the Rayleigh quotient
\begin{align}
    \varepsilon_i := \sqrt{ \frac{Q(c_i,\widetilde{H}_i)}{Q(c_i,H)} },
    \label{eqn:epsdef}
\end{align}
where $\widetilde{H}_i(\lambda) = H(\lambda) (\lambda - \widetilde{\lambda}_i)^2$, for the same $c_i$ computed above. Then $\varepsilon_{i}^{2}$ is a weighted average of $(\lambda - \widetilde{\lambda}_i)^2$ and hence there exists a cuspidal eigenvalue $\lambda \in [\widetilde{\lambda}_i - \varepsilon_i, \widetilde{\lambda}_i + \varepsilon_i]$.

Next we prove completeness of the eigenvalues, i.e. prove that we have not missed any. We choose a test function $H^*(\lambda)$ that is positive and monotonically decreasing for $\lambda > 0$. Then $H^*(\lambda) \geq H^*(\widetilde{\lambda}_i + \varepsilon_i)$ for all $\lambda \in [\widetilde{\lambda}_i - \varepsilon_i, \widetilde{\lambda}_i + \varepsilon_i]$. Hence any eigenvalue $\lambda$ that is not contained in $\bigcup_i [\widetilde{\lambda}_i - \varepsilon_i, \widetilde{\lambda}_i + \varepsilon_i]$ must satisfy
\begin{align*}
    H^*(\lambda) \leq t(1,H^*) + \sum_i H^*(\widetilde{\lambda}_i + \varepsilon_i).
\end{align*}
Here the second sum ranges over all $i$ such that $[\widetilde{\lambda}_i - \varepsilon_i, \widetilde{\lambda}_i + \varepsilon_i]$ does not overlap the corresponding interval for any smaller value of $i$. Since $H^*$ is monotonic, this determines numbers $\delta_i > 0$ such that $|\lambda_i - \widetilde{\lambda}_i| \leq \varepsilon_i$ and $|\lambda_j - \widetilde{\lambda}_i|\geq \delta_i$ for $j \in \NN \backslash \lbrace i \rbrace$. Note that this approach only works well if the $\lambda_i$ turn out to be distinct and well separated. It is conjectured that the Laplacian spectrum is simple for squarefree level and trivial character, with Poissonian spacing statistics. There exists some theoretical and numerical evidence for this, namely from \cite{LuoSarnakspacing} and \cite{Steil} respectively. For this algorithm, we will see from the data that this will be the case.

Finally, we consider the Hecke eigenvalues. For $j \geq 1$ and any sequence $\lbrace c(m) \rbrace_{m = 1}^{M}$, define
\begin{align*}
    L_j(c) = \sum_{m=1}^{M} c(m) a_j(m).
\end{align*}
Let $H, \widetilde{H}_i$ be as above. Then
\begin{align*}
    \left( \sum_{j \neq i} L_j(c_i) a_j(n) H(\lambda_j) \right)^2 &\leq \sum_{j \neq i} L_j(c_i)^2 H(\lambda_j) \sum_{j=1}^{\infty} (a_j(n))^2 H(\lambda_j)\\
    &\leq \delta_i^{-2} Q(c_i, \widetilde{H}_i) Q(e_n, H) = \varepsilon_i^2 \delta_i^{-2} Q(c_i, H) Q(e_n, H), 
\end{align*}
where $e_n(m) = 1$ if $m=n$ and $0$ otherwise. Thus, defining
\begin{align}
    \eta_{i,n} = \frac{\varepsilon_i}{\delta_i} \sqrt{Q(c_i, H) Q(e_n,H)} \quad \textup{ and } \quad W_i = L_i(c_i) H(\lambda_i),
    \label{eqn:etadef}
\end{align}
we have
\begin{align*}
    A_i(n) := a_i(n) W_i = \sum_{m=1}^{M} c_i(m) \sum_{d | (m,n)} t \left( \frac{mn}{d^2}, H\right) + \beta_{i,n} \eta_{i,n},
\end{align*}
where $\beta_{i,n}$ is some real constant that depends on $i$ and $n$ and satisfies $|\beta_{i,n}| \leq 1$. We can use this to compute $a_i(n)$, with $(n,N) = 1$, by using the fact that $a_i(1)=1$ to compute $W_i$ to a proven accuracy.

In practice, we will choose one test function $H$ that is both positive and monotonically decreasing and use this throughout.

\subsection{Computing $a_n$ for $(n,N) > 1$ for squarefree level $N$}
\label{subsec:hecke_eigenvals_coprime}
Let $f$ be a primitive Maass newform of squarefree level $N$, Laplace eigenvalue $\lambda = 1/4 + R^2$ and trivial character, with Fourier coefficients $a_n$. By Atkin--Lehner theory \cite{AtkinLehner} (see \cite[\S 1.2.6]{stromberg} for non-holomorphic case), for each prime $p \mid N$ we have $a_p = \pm 1/\sqrt{p}$. Moreover, defining $w = \mu(N) \sqrt{N} \prod_{p \mid N} a_p = \prod_{p|N} \sign(-a_p)$, we have $f(z) = w f(-1/Nz)$. Hence, we just need to find the signs of the $a_p$ for $p \mid N$, and then use the Hecke relations to find all $a_n$ for $(n,N)>1$.

Suppose first that $f$ is even, so its Fourier expansion is of the form
\begin{align*}
    f(z) = \sum_{n=1}^{\infty} \frac{a_n}{\sqrt{n}} W_{iR}(2 \pi n y) \cos(2 \pi n x),
\end{align*}
where $W_s(y) := \sqrt{y} K_s(y)$ and $K_s(y)$ is the K-Bessel function. Substituting $z=iy$ into the relation $f(z) = w f(-1/Nz)$, we have
\begin{align}
    \sum_{n=1}^{\infty} \frac{a_n}{\sqrt{n}} \left( W_{iR}(2 \pi n y) - w W_{iR} \left( \frac{2 \pi n}{Ny}\right) \right) = 0.
    \label{eqn: Weqneven}
\end{align}
If $w = -1$ then taking $y = 1/\sqrt{N}$ in \eqref{eqn: Weqneven} yields
\begin{align*}
    \sum_{n=1}^{\infty} \frac{a_n}{\sqrt{n}} W_{iR} \left( \frac{2 \pi n}{\sqrt{N}} \right) = 0.
\end{align*}
If $w = 1$ then taking $y = \sqrt{2/N}$ in \eqref{eqn: Weqneven} yields
\begin{align*}
    \sum_{n=1}^{\infty} \frac{a_n}{\sqrt{n}} \left( W_{iR} \left( \frac{2 \pi n \sqrt{2}}{\sqrt{N}} \right) - W_{iR} \left( \frac{\pi n \sqrt{2}}{\sqrt{N}} \right)\right)= 0.
\end{align*}
Now suppose $f$ is odd, so its Fourier expansion takes the form
\begin{align*}
    f(z) = \sum_{n=1}^{\infty} \frac{a_n}{\sqrt{n}} W_{iR}(2 \pi n y) \sin(2 \pi n x).
\end{align*}
In this case plugging in $z=iy$ would only give the trivial relation $0=0$, so instead we first differentiate with respect to $x$. For this we consider
\begin{align*}
    \frac{\partial}{\partial x} (f(z) - w f(-1/Nz))|_{z=iy} = 0.
\end{align*}
After some computation this yields
\begin{align}
    \sum_{n=1}^{\infty} a_n \sqrt{n} \left(W_{iR}(2 \pi n y) + \frac{w}{N y^2} W_{iR} \left( \frac{2 \pi n}{Ny} \right) \right)=0.
    \label{eqn: Weqnodd}
\end{align}
If $w=1$ then taking $y=1/\sqrt{N}$ in \eqref{eqn: Weqnodd} yields
\begin{align*}
    \sum_{n=1}^{\infty} a_n \sqrt{n} W_{iR} \left( \frac{2 \pi n}{\sqrt{N}} \right) = 0.
\end{align*}
If $w=-1$ then taking $y=\sqrt{2/N}$ in \eqref{eqn: Weqnodd} yields
\begin{align*}
    \sum_{n=1}^{\infty} a_n \sqrt{n} \left( W_{iR} \left( \frac{2 \pi n \sqrt{2}}{\sqrt{N}} \right) - \frac{1}{2} W_{iR} \left( \frac{\pi n \sqrt{2}}{\sqrt{N}} \right) \right) = 0.
\end{align*}
In summary, if we define
\begin{align*}
    W(y) = \begin{dcases}
        W_{iR}(y) & \text{if $f$ is even and $w=-1$,} \\
        W_{iR}(y\sqrt{2}) - W_{iR}(y/\sqrt{2}) & \text{if $f$ is even and $w=1$,}\\
        \frac{y \sqrt{N}}{2 \pi} W_{iR}(y) & \text{if $f$ is odd and $w=1$,}\\
        \frac{y \sqrt{N}}{2 \pi} \left(W_{iR}(y\sqrt{2}) - \frac{1}{2}W_{iR}(y/\sqrt{2})\right) & \text{if $f$ is odd and $w=-1$,}
    \end{dcases}
\end{align*}
then
\begin{align*}
    \sum_{n=1}^{\infty} \frac{a_n}{\sqrt{n}} W \left( \frac{2 \pi n}{\sqrt{N}} \right) = 0.
\end{align*}

Now computationally we will only have accurate approximations of $a_n$ for $n \leq M$, so we must truncate the above sums at $M$ and estimate the error incurred.
Using the current best estimate towards to Ramanujan--Petersson conjecture from Kim--Sarnak \cite{kimsarnak}, we have $|a_p| \leq p^{7/64} + p^{-7/64}$, which implies $|a_n/\sqrt{n}| \leq \theta \approx 1.758$. We also have that
\begin{align*}
    |W_{iR}(y)| \leq  \sqrt{\frac{\pi}{2}} e^{-y} \text{ for } y > 0.
\end{align*}
With both of these results we can easily find bounds for the tails of the sums and obtain
\begin{align*}
    \left|  \sum_{n=M+1}^{\infty} \frac{a_n}{\sqrt{n}} W \left( \frac{2 \pi n}{\sqrt{N}} \right) \right|\leq \begin{dcases}
        \theta \sqrt{\frac{\pi}{2}} \frac{\exp\left(- \frac{2 \pi M}{\sqrt{N}} \right)}{\exp\left(\frac{2 \pi }{\sqrt{N}} \right)-1} & \text{if $f$ is even and $w=-1$,} \\
        2\theta \sqrt{\frac{\pi}{2}} \frac{\exp\left( -\frac{\pi M \sqrt{2}}{\sqrt{N}} \right)}{\exp\left(\frac{\pi \sqrt{2}}{\sqrt{N}} \right)-1} & \text{if $f$ is even and $w=1$,}\\
        \theta \sqrt{\frac{\pi}{2}} \frac{\left((M+1) \exp \left( \frac{2 \pi }{\sqrt{N}} \right) - M\right)}{ \exp \left(\frac{2 \pi M}{\sqrt{N}} \right) \left(\exp \left( \frac{2 \pi }{\sqrt{N}} \right) - 1\right)^2} & \text{if $f$ is odd and $w=1$,}\\
        \frac{3 \theta}{2} \sqrt{\frac{\pi}{2}} \frac{ \left((M+1) \exp \left( \frac{\pi \sqrt{2}}{\sqrt{N}} \right) - M\right)}{\exp \left(\frac{\pi M \sqrt{2}}{\sqrt{N}} \right)\left(\exp \left( \frac{\pi \sqrt{2}}{\sqrt{N}} \right) - 1\right)^2} & \text{if $f$ is odd and $w=-1$.}
    \end{dcases}
\end{align*}
To find the signs of the $a_p$ for $p \mid N$ we just test every combination of $\pm 1$ for the signs of the $a_p$, then use this to compute $w$ and the corresponding sum from the above cases. Heuristically, we expect only one of these sums to be within the error derived. When there is only one sum within the errors, we can say that the result is rigorous. We then take the signs of the $a_p$ for $p \mid N$ and $w$ from that sum. In practice we see this works well, provided the Laplace eigenvalue and Hecke eigenvalues are computed to a high enough precision.
\section{The Selberg Trace Formula for squarefree level $N$}
\label{sec:traceformula}
In the algorithm given in Section \ref{sec:algorithm}, an essential tool we need is an explicit version of the Selberg trace formula with Hecke operators. Currently, this has only been derived for squarefree level by Str\"ombergsson in \cite{Andreaspre}. For our computation, we rewrite this in the following form, following the steps of Proposition 2.2 in \cite{AndyMin}.

\begin{theorem}[The Selberg trace formula for Maass newforms for squarefree level and trivial character]
Fix $\delta > 0$, let $h(t)$ be a even analytic function on the strip $\lbrace t \in \CC :\Imag(t) \leq \frac{1}{2} + \delta \rbrace$ such that $h(r) \in \RR$ for $r \in \RR$ and $h(r) = O((1 + |r|^2)^{-1-\delta})$. Define $g$ as the Fourier transform of $h$ given by
\begin{align*}
	g(u) = \frac{1}{2 \pi} \int_{-\infty}^{\infty} h(r) e^{-iru} dr.
\end{align*}

Let $\lbrace f_j \rbrace$ be a sequence of normalised Hecke eigenforms of squarefree level $N$, with Laplacian eigenvalues $\lambda_j = \frac{1}{4} + r_j^2$ and respective Hecke eigenvalues $a_j(n)$.

Then, for $(N,n)=1$ we have
\begin{align*}
	&\frac{\mu(N) \sigma_1(|n|)}{\sqrt{|n|}} h \left( \frac{i}{2} \right) + \sum_{j > 0} h(r_j) a_j(n)\\
	&= \underset{\substack{t \in \ZZ\\\sqrt{D} = \sqrt{t^2 - 4n} \not\in \QQ}}{\sum} c_N(D) \cdot \begin{dcases}
    g \left( \log \left( \frac{(|t| + \sqrt{D})^2}{4 |n|} \right)\right) & \text{if } D > 0,\\
    \frac{\sqrt{|D/4n|}}{2 \pi} \int_{-\infty}^{\infty} \frac{g(u) \cosh(u/2)}{\sinh^2(u/2) + |D/4n|}du & \text{if } D < 0  \end{dcases}\\
    &+\Lambda(N) \underset{\substack{ad = n \\ a>0 \\ a \neq d}}{\sum} \frac{ g\left( \log \left| \frac{a}{d} \right|\right)}{(N^{\infty},|a-d|)} - 2\Lambda(N) \underset{\substack{ad = n \\ a>0}}{\sum} \sum_{r=0}^{\infty} N^{-r} g\left( \log \left| \frac{a}{d} \right| - 2r \log(N)\right) \\
    &+\begin{dcases}
   -\frac{\prod_{p|N}(p-1)}{12 \sqrt{n}} \int_{-\infty}^{\infty} \frac{g'(u)}{\sinh\left( \frac{u}{2} \right)}du & \textup{ if } \sqrt{n} \in \ZZ, \\
   0 & \textup{ otherwise,}
   \end{dcases}\\
\end{align*}
where
\begin{align*}
    c_N(D) = L(1,\psi_D) \prod_{p|N} (\psi_d(p) - 1) = \frac{L(1, \psi_d)}{l} \prod_{p|N} (\psi_d(p) - 1) \prod_{p | l} \left[ 1 + (p - \psi_d(p)) \frac{(l,p^{\infty}) - 1}{p-1} \right],
\end{align*}
with $D = dl^2$, $l>0$, $d$ a fundamental discriminant and $\psi_d(p) = \left( \frac{d}{p} \right)$. Here $(l,p^{\infty})$ denotes the largest power of $p$ that divides $l$.
\end{theorem}

\begin{remark}
	We refer to the terms in the sum with $D>0$ as the hyperbolic terms and the terms $D<0$ as the elliptic terms. The terms that are multiplied by the von Mangoldt function $\Lambda(N)$ we call the parabolic terms, and the term when $\sqrt{n} \in \ZZ$ we call the identity term.
\end{remark}

\subsection{Computational remarks}
The main numerical bottleneck of computing the trace formula is from the contribution of the hyperbolic terms, which involves computing the class number and regulator of $\QQ(\sqrt{D})$. For numerical stability, it is best to consider a test function $g$ that is compactly supported. This allows one to compute the terms on the geometric side to arbitrary precision with a fixed finite list of class numbers. Precisely we would need class numbers $h_{\QQ(\sqrt{D})}$ for $D = t^2 - 4n < (2n\cosh(X/2))^2$. For our computation, we used Pari \cite{PARI} to compute these real class numbers and regulators, and verified the calculations with \cite{Jacobsonclassnumbers}.

We can also get a bonus increase in the precision by considering the parity of the forms separately. The traces are given by $\frac{1}{2} (t(n,h) + t(-n,h))$ and $\frac{1}{2} (t(n,h) - t(-n,h))$ for the even and odd forms respectively.

Computing the integrals appearing in the elliptic terms to arbitrary precision can also be challenging given the large number of them appearing for values of $D$ and $n$. We can remedy this by noting that $|D/4n| \in (0,1]$ for $D = t^2 - 4n < 0$, and considering the integrals as functions $f:(0,1] \to \RR$ defined by
\begin{align*}
	f(x) = \int_{0}^{\infty} \frac{g(u) \cosh(u/2)}{\sinh^2(u/2) + x}du.
\end{align*}
This function is analytic with respect to the variable $x$, hence we can approximate this integral with a Taylor series, where the only integrals we need to compute are given in the Taylor coefficients. Explicitly, for $x$ near $x_0$, we can approximate $f(x)$ by
\begin{align*}
    f(x) = \sum_{k=0}^{K} \frac{f^{(k)}(x_0)}{k!} (x-x_0)^k + R_K(x),
\end{align*}
where $R_K(x)$ is the error term given by
\begin{align*}
    R_K(x) = \frac{f^{(K+1)}(\xi)}{(K+1)!} (x-x_0)^{K+1},
\end{align*}
for some $\xi$ in the closed interval between $x$ and $x_0$. To find the Taylor coefficients, we use Leibniz's integral rule to get
\begin{align*}
    \frac{d^k}{dx^k} f(x) = k! (-1)^k  \int_0^{\infty} \frac{g(u) \cosh(u/2)}{(\sinh^2(u/2) + x)^{k+1}}du.
\end{align*}
To bound the error term, let $\xi \in [x_0,x]$ and $M_g = \max_{y \in [0,\infty)} |g(y)|$. Then
\begin{align*}
    |f^{(K+1)}(\xi)| &= (K+1)! \left| \int_0^{\infty} \frac{g(u) \cosh(u/2)}{(\sinh^2(u/2) + \xi)^{K+2}}du \right|
    \leq M_g (K+1)! \int_0^{\infty} \frac{\cosh(u/2)}{(\sinh^2(u/2) + \xi)^{K+2}}du.
\end{align*}
Here we have that
\begin{align*}
    \int_0^{\infty} \frac{\cosh(u/2)}{(\sinh^2(u/2) + \xi)^{K+2}}du = \pi \xi^{-3/2-K} \prod_{k=1}^{K+1} \left( \frac{2k-1}{2k}\right).
\end{align*}
Hence we can bound the error term in the Taylor series by
\begin{align*}
    |R_K(x)| \leq \frac{\pi M_g}{\sqrt{x_0}} \left|1 - \frac{x}{x_0}\right|^{K+1} \prod_{k=1}^{K+1} \left( \frac{2k-1}{2k}\right).
\end{align*}
To compute all the elliptic integrals, we shall need to choose the sample points for our Taylor series, such that it minimises the number of Taylor coefficients that are needed to be computed. Since there is a singularity at $x=0$, it is best for us to choose our sampling points geometrically, that is $x_j = c^{-j}$ for some $c>1$. Suppose, we take $K$ terms of a Taylor expansion around the point $x_j$, we can see that error is of size about $|1-x/x_0|^{K}$. For our sample points, we have
\begin{align*}
    \left|1 - \frac{x}{x_0}\right| \leq \left( \frac{c-1}{c+1} \right),
\end{align*}
hence the worst our error could be is $\left( \frac{c-1}{c+1} \right)^K$. Note, that given $x$ we can choose $j = \lceil \log_c(\frac{2}{(c+1)x}) \rceil$. Thus to choose the number of sampling points needed, we just consider the smallest value of $x$ that we could feasibly have.

We see that the number of sample points is about $\log_c n$, where $n$ is the largest Hecke operator we shall need to consider. So in total we have to compute about $K \log_c n$ integrals, and we want to minimise this with respect to the constraint that $\left( \frac{c-1}{c+1} \right)^K < \varepsilon$ for some fixed error tolerance $\varepsilon$. This surprisingly has the exact solution with $c = 1+\sqrt{2}$ and $K = \log_c(1/\varepsilon)$.

\section{Choice of test function}
\label{sec:testfunc}
As stated in Sections \ref{sec:algorithm} and \ref{sec:traceformula}, we will want a test function that is even, positive and monotonically decreasing. Moreover, to aid in computations, we will also want $g$, the Fourier transform of $h$, to be compactly supported. This will make all the integrals and sums on the geometric side have finite bounds which will help when implementing the algorithm.

\subsection{Candidate test function}
A good initial function to consider is powers of the $\sinc(x) = \sin(x)/x$ function. For even powers, this is a positive even function with a compactly supported Fourier transform. However, this function is not monotonically decreasing. To remedy this we consider the test function
\begin{align*}
    h_1(t) = \frac{\pi^2}{\pi^2 + 4} \left[ \sinc^2\left( \frac{t}{2} \right) + \frac{1}{2} \sinc^2\left( \frac{t-\pi}{2} \right) + \frac{1}{2} \sinc^2\left( \frac{t+\pi}{2} \right) \right],
\end{align*}
and let $h_d(t) = h_1(t)^d$ for $d \in \NN$. Then $h_d(t)$ is a positive, even and monotonically decreasing function on $\RR_{>0}$, satisfying $h_d(0)=1$ and
\begin{align*}
    h_d(t) \sim \left( \frac{4 \pi^2}{\pi^2 + 4} \right)^d t^{-2d},
\end{align*}
as $|t| \to \infty$. Moreover, its Fourier transform
\begin{align*}
	g_d(x) = \frac{1}{\pi} \int_0^{\infty} h_d(t) \cos(t x) dt,
\end{align*}
is compactly supported on $[-d,d]$. For a fixed $d$ we can express $g_d$ in the form
\begin{align*}
    g_d(x) = \sum_{m \in \lbrace -1,0,1 \rbrace} A_m(x) e^{\pi i m x},
\end{align*}
where
\begin{align*}
    A_m(x) = A_{m,j} \left(x - j - \frac{1}{2}\right) \quad \textup{ for } x \in [j,j+1), j \in \lbrace -d, \ldots, d-1 \rbrace,
\end{align*}
for some $A_{m,j} \in \CC[x]$ satisfying $A_{m,-1-j}(x) = A_{-m,j}(-x) = \overline{A_{m,j}(-x)}$. Note that all the $A_{m,j}$ are determined by those with $m \in \lbrace 0,1 \rbrace$ and $j \in \lbrace 0,\ldots,d-1 \rbrace$.

Specifically, for $d=1$, we have
\begin{align*}
    A_{0,0}(x) = \frac{\pi^2}{\pi^2 + 4} \left(\frac{1}{2}-x\right) \quad \textup{ and  } \quad A_{1,0}(x) = \frac{1}{2} A_{0,0}(x).
\end{align*}

For $d>1$, we compute the functions using convolutions. More explicitly, suppose we are given functions
\begin{align*}
    A(x) = \sum_{m \in \{ -1,0,1 \}} A_m(x) e^{\pi i m x} \quad \textup{ and  } \quad B(x) = \sum_{m \in \{ -1,0,1 \}} B_m(x) e^{\pi i m x},
\end{align*}
and we wish to compute their convolution $C = A*B$, which is again a function of the same form. For a set $S$, we define the indicator function $\textbf{1}_S(x) = 1$ if $x \in S$ and $0$ if $x \not\in S$. It suffices to consider the constituent functions
\begin{align*}
    A_{m,j}\left( x-j-\frac{1}{2} \right) e^{\pi i m x}\textbf{1}_{[j,j+1)}(x) \quad \textup{ and } \quad B_{n,k}\left( x-k-\frac{1}{2} \right) e^{\pi i n x}\textbf{1}_{[k,k+1)}(x),
\end{align*}
with convolution
\begin{align*}
    \int_{\RR} A_{m,j}\left( y-j-\frac{1}{2} \right) e^{\pi i m y}\textbf{1}_{[j,j+1)}(y) B_{n,k}\left( x-y-k-\frac{1}{2} \right) e^{\pi i n (x-y)}\textbf{1}_{[k,k+1)}(x-y)dy.
\end{align*}
Consider $x \in [j+k+\delta,j+k+\delta+1)$ for some $\delta \in \{ 0,1 \}$, and let $t = x - \left( j+k+\delta + \frac{1}{2} \right)$. We make the change of variable $y \mapsto y + j + \frac{1}{2}$ to get
\begin{align}
    &\int_{\RR} A_{m,j}(y) e^{\pi i m \left(y+j+\frac{1}{2}\right)} \textbf{1}_{\left[-\frac{1}{2}\right.,\left.\frac{1}{2}\right)}(y)B_{n,k} \left( t + \delta - \frac{1}{2} -y \right) e^{\pi i n \left(x - y - j - \frac{1}{2}\right)} \textbf{1}_{(t+\delta-1,t+\delta]}(y)dy \nonumber \\
    &= e^{\pi i (m-n) \left( j + \frac{1}{2} \right) + \pi i n x} (-1)^{\delta} \int_{\delta - \frac{1}{2}}^{t} A_{m,j}(y) B_{n,k} \left( t + \delta - \frac{1}{2} - y\right) e^{\pi i (m-n)y} dy.
    \label{eqn:convolutioneqn}
\end{align}
When $m \neq n$ we apply repeated integration by parts to see that \eqref{eqn:convolutioneqn} becomes
\begin{align*}
    &e^{\pi i (m-n) \left( j + \frac{1}{2} \right) + \pi i n x} (-1)^{\delta} \sum_{r=0}^{\deg{A_{m,j}}} \sum_{s=0}^{\deg{B_{n,k}}} \frac{(-1)^s \binom{r+s}{s}}{(-\pi i (m-n))^{r+s+1}}\\
    &\cdot \left( A_{m,j}^{(r)} \left( \delta - \frac{1}{2} \right) B_{n,k}^{(s)}(t) e^{\pi i (m-n) \left( \delta - \frac{1}{2} \right)} - A_{m,j}^{(r)}(t) B_{n,k}^{(s)} \left( \delta - \frac{1}{2} \right) e^{\pi i (m-n)t} \right)\\
    &= (-1)^{(m-n+1) \delta} \sum_{r=0}^{\deg{A_{m,j}}} \sum_{s=0}^{\deg{B_{n,k}}} \frac{(-1)^s \binom{r+s}{s}}{(-\pi i (m-n))^{r+s+1}}\\
    &\cdot \left( A_{m,j}^{(r)} \left( \delta - \frac{1}{2} \right) B_{n,k}^{(s)}(t) (-1)^{(m-n)j} e^{\pi i n x} - A_{m,j}^{(r)}(t) B_{n,k}^{(s)} \left( \delta - \frac{1}{2} \right) (-1)^{(m-n)k}e^{\pi i m x} \right).
\end{align*}
Note that this will contribute to both the $C_{m,j+k+\delta}$ and $C_{n,j+k+\delta}$ terms.

When $m=n$, we define polynomials $P_{\delta,l} \in \CC[y]$ such that $P_{\delta,0} = A_{m,j}(y)$ and
\begin{align*}
    P_{\delta,l} = \int_{\delta - \frac{1}{2}}^{y} P_{\delta,l-1}(u)du,
\end{align*}
for $l \geq 1$. Then applying integration by parts, \eqref{eqn:convolutioneqn} becomes
\begin{align*}
    (-1)^{\delta} \sum_{l=1}^{\deg{B_{n,k}}+1} B_{n,k}^{(l-1)} \left( \delta - \frac{1}{2} \right) P_{\delta,l}(t) e^{\pi i m x}.
\end{align*}
\subsection{Optimising the test function}
\label{subsec:optimisetestfunc}
We wish to optimise the decay of the test function for certain given constants such that we maximise the precision with which we compute the trace formula. Suppose we aim for a final precision of $B$ bits. Due to the square roots in \eqref{eqn:epsdef} and \eqref{eqn:etadef}, we must consider terms larger than $2^{-2B}$ to be significant, and use a working precision of at least $2B$ bits. Let $X \in \RR_{>0}, d \in \NN$ and consider the test function
\begin{align}
	h(r) = h_d \left( \frac{Xr}{d} \right).
	\label{eqn:testfuncdil}
\end{align}
From this we see that $g$, the Fourier transform of $h$, is compactly supported on $[-X,X]$. We take the edge of the precision window to be the point $R_{\rm max}$ at which 
\begin{align}
    h(R_{\rm max}) = h_1(XR_{\rm max}/d)^d = 2^{-2B}.
    \label{eqn: h_1(XR/d)}
\end{align}

Fix a level $N$. Let $M$ be the number of level $N$ newforms with trivial character, fixed parity and Laplace eigenvalue $\lambda \leq \frac{1}{4} +R_{\rm max}^2$ and let $D_{\rm max}$ be the largest size of discriminant appearing in the hyperbolic sum. The value $M$ will control the size of the matrices appearing in the linear algebra and $D_{\max}$ will control how many hyperbolic terms will appear. We want the ability to choose these values since these are the main sections of the algorithm that are constrained by external factors. For example, we will only have a list of class numbers up to a certain limit that that could feasibly be computed. The idea of this section is to first fix $N, M$ and $D_{\rm max}$, then find $R_{\rm max},X$ and $d$ such that it maximises the precision $B$.

So fix $N,M$ and $D_{\rm max}$. To find $R_{\rm max}$, we have from \cite{Risager} that
\begin{align*}
	M = \frac{R_{\rm max}^2}{24}N + O( \sqrt{\lambda} \log{\lambda}),
\end{align*}
which we can rearrange to compute $R_{\rm max}$ by
\begin{align*}
	R_{\rm max} \approx \sqrt{\frac{24M}{N}}.
\end{align*}
To find $X$, we use the fact that $g$ is compactly supported on $[-X,X]$ and hence, we have that
\begin{align*}
	D_{\rm max} = \left(2M\cosh \left( \frac{X}{2} \right)\right)^2,
\end{align*}
which we can rearrange to compute $X$ by
\begin{align*}
	X = 2\cosh^{-1} \left( \frac{\sqrt{D_{\rm max}}}{2M} \right).
\end{align*}
Once we have values for $R_{\rm max}$ and $X$, we can find $d$ by first rearranging \eqref{eqn: h_1(XR/d)} to obtain
\begin{align*}
-\log_2 \left( h_1 \left( \frac{XR_{\rm max}}{d} \right) \right)d = 2B.
\end{align*}
We can now find a $d$ which maximises the left side of this equation, which in turn will maximise our final precision $B$. Note that since $d \in \NN$, we can find the maximum by sampling the left side of the equation over integer values of $d$ and choosing the largest value.

Thus, once we have computed these values, the test function we use for the computation is given by \eqref{eqn:testfuncdil}. In practice, when choosing the level $N$, we pick $N$ to be the largest level we are computing with and use this test function for all smaller levels as well.

\section{Computational results}
\label{sec:compres}

\subsection{Computing the forms}
We implemented this algorithm in the C programming language, predominately using the ball-arithmetic library \textbf{Arb} \cite{arb} throughout our computations to manage round-off errors. For the main computation, following the notation from Section \ref{subsec:optimisetestfunc}, we chose the numbers $D_{\rm max} = 10^9,M = 2000$ and the maximum level we consider is $N = 105$. Using SageMath \cite{sagemath}, we find $X \approx 5.51341, R_{\rm max} \approx 21.38089, d=13$ and $2B \approx 63$.

With these numbers, we computed a total of $33214$ Laplace eigenvalues of Maass cusp forms, each with all Hecke eigenvalues $a_n$ with $n \leq 2000$ and $(n,N)=1$, for squarefree levels $2 \leq N \leq 105$. The range of the $\varepsilon_i$'s computed is between $10^{-15}$ and $10^{-2}$. Of these forms $17243$ are even and $15971$ are odd.

Of these Laplace eigenvalues, we proved completeness for $16207$ of them and hence, their Hecke eigenvalues have rigorous error bounds. We could only compute completeness for all prime levels $2 \leq N \leq 67$ and all composite squarefree levels $6 \leq N \leq 105$ due to the precision of the computed trace formula values in the linear algebra. Each of these complete Laplace eigenvalues will correspond to a provably unique Maass cusp form. Of these forms $8419$ are even and $7788$ are odd.

We observed that the closest distance between two Maass forms in the completed range was approximately $3 \times 10^{-6}$ from the level $23$ Laplace eigenvalues of $10.85166055 \ldots$ and $10.8516021 \ldots$. The closest distance between two even forms was approximately $1.4 \times 10^{-5}$ from the level $53$ Laplace eigenvalues of $5.876312 \ldots$ and $5.876299 \ldots$. The closest distance between two odd forms was approximately $3 \times 10^{-6}$ from the level $55$ Laplace eigenvalues of $8.350572 \ldots$ and $8.350569 \ldots$.

The entire computation took just under two weeks of time on 64 cores of 2.5GHz AMD Opteron processors. As predicted, the computation was dominated by computing the hyperbolic terms. 
\subsection{Ramanujan--Petersson conjecture}
The Ramanujan--Petersson conjecture states that for prime $p$, the $p$th Fourier coefficient $a_p$ for a Maass cusp form on $\Gamma_0(N)$ should satisfy $|a_p| \leq 2$. For the data we computed, we verified this was true for all Hecke eigenvalues with $p \leq 2000$ for $13271$ of our Maass forms that we proved completeness for. 

\subsection{Sato--Tate conjecture}
\begin{figure}[ht]
    \centering
    \begin{subfigure}[t]{0.485\textwidth}
        \centering
        \def\svgwidth{\linewidth}
        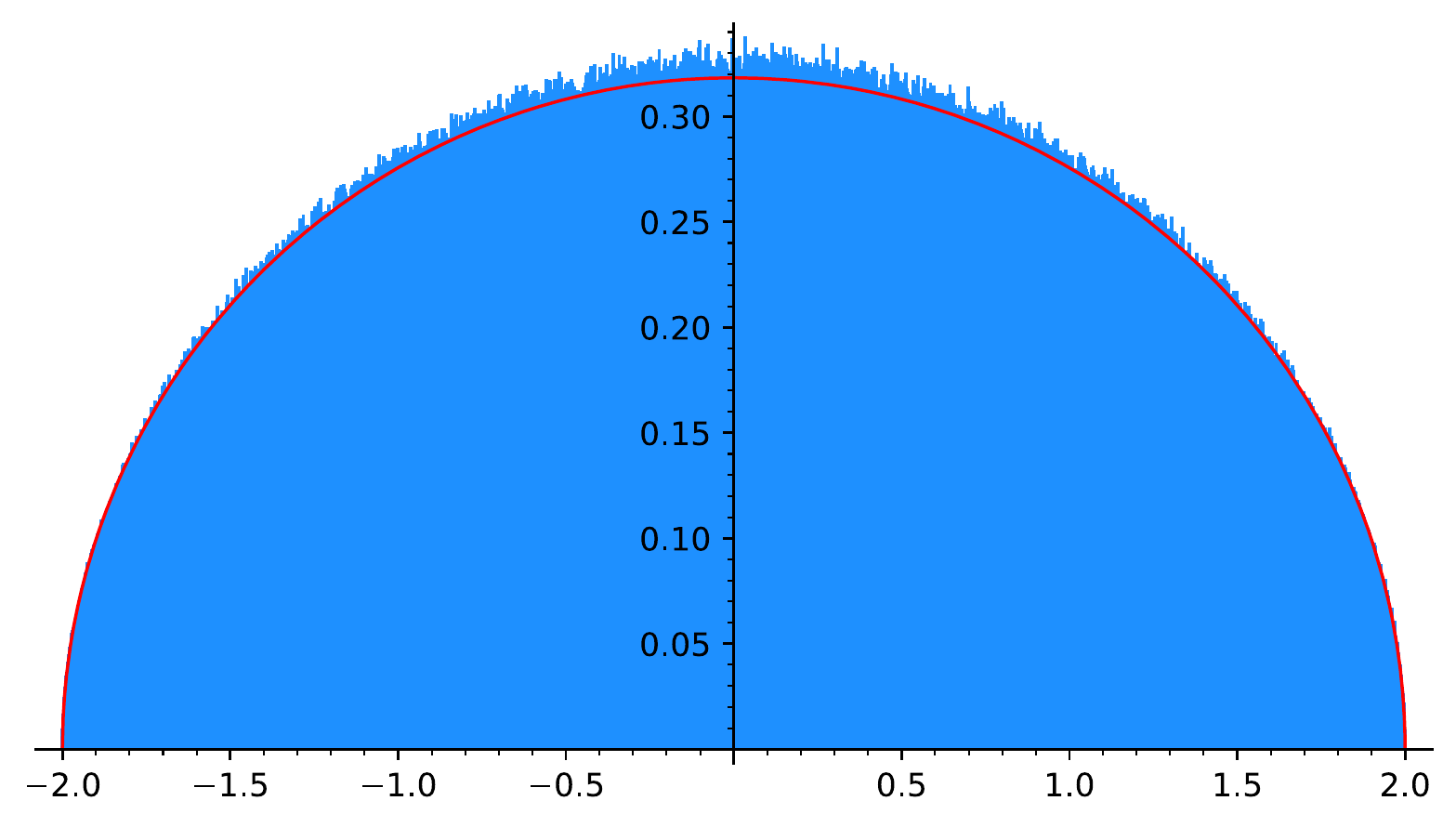
        \label{fig: satotatep}
    \end{subfigure}
    ~
    \begin{subfigure}[t]{0.485\textwidth}
        \centering
        \def\svgwidth{\linewidth}
        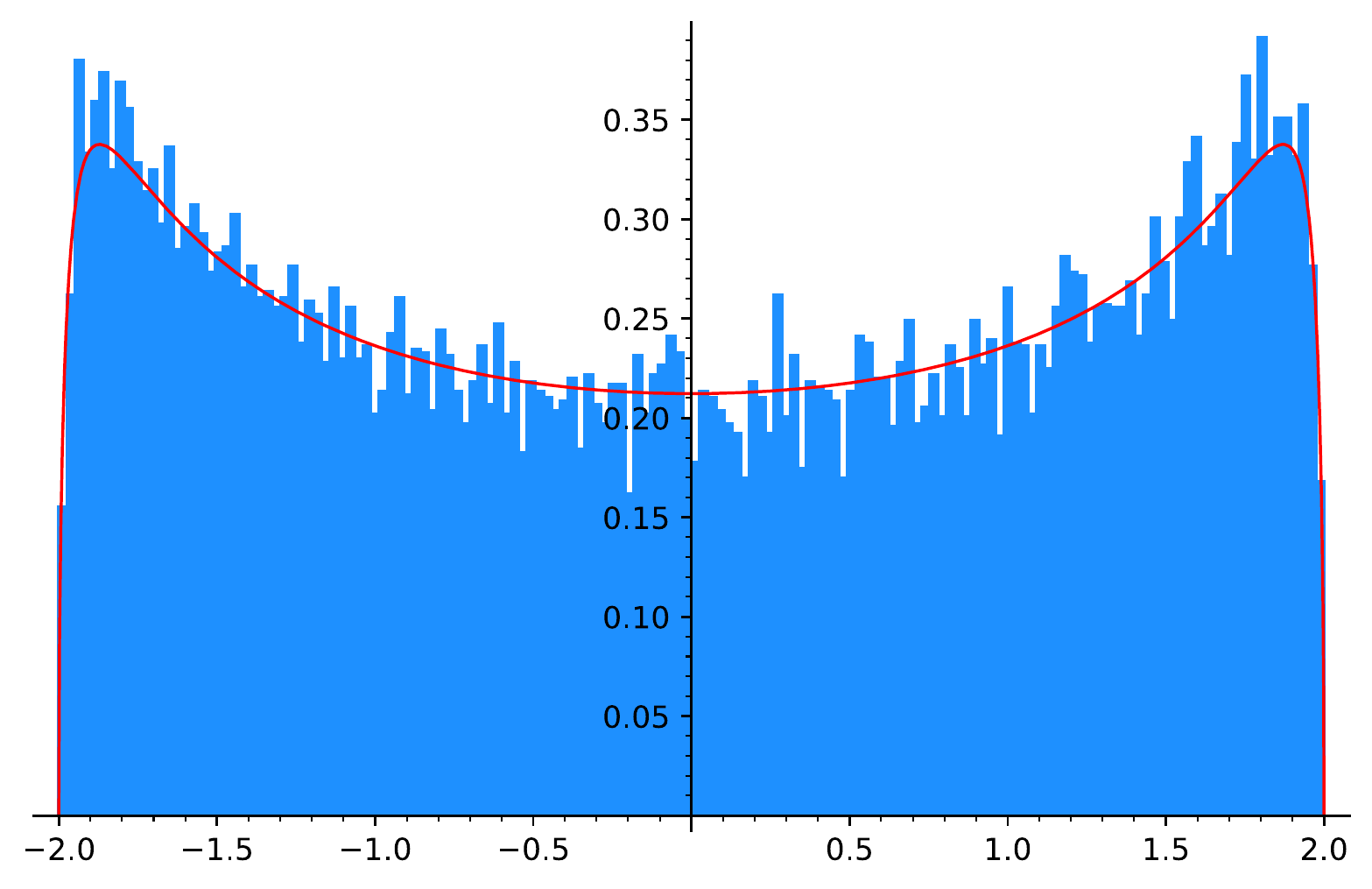
        \label{fig: satotat2}
    \end{subfigure}
    \caption{Comparison of our data to predicted distributions. The left-hand figure concerns the distribution of the classical Sato--Tate conjecture; the histogram has $10003411$ data points in $3162$ bins. The right-hand figure compares our data to Sarnak's theorem \cite{sarnak} for $a_2$; the histogram has $23806$ data points in $154$ bins.}
    \label{fig: histos}
\end{figure}

The Sato--Tate conjecture is a statistical conjecture about the asymptotic distribution of eigenvalues $a_p$ of Hecke operators $T_p$ for primes $p$. It states that the $a_p$ should be asymptotically distributed with the Sato--Tate measure given by
\begin{align*}
	\mu_{\infty} = \frac{1}{\pi} \sqrt{1 - \frac{x^2}{4}} dx,
\end{align*}
as $p \to \infty$. A related result, proven by Sarnak in \cite{sarnak}, states that instead if we fix a prime $p \nmid N$ and let the level tend to infinity, then the points $a_p$ of these forms are asymptotically distributed by the measure
\begin{align*}
	\mu_p = f_p \mu_{\infty},
\end{align*}
where
\begin{align*}
	f_p(x) = \frac{p+1}{(p^{1/2} + p^{-1/2})^2 - x^2},
\end{align*}
for $x \in [-2,2]$. As an example, for $p=2$, the points should be distributed asymptotically with respect to
\begin{align}
	\mu_2 = \frac{3 \sqrt{4-x^2}}{9 - 2x^2} \frac{dx}{\pi}.
	\label{eqn: mu2}
\end{align}
We used the Maass form data to create Figure \ref{fig: histos}, which illustrates a strong connection to the predicted result of the Sato--Tate conjecture and the result proven by Sarnak.

\subsection{$L$-function and the Riemann hypothesis}
Let $f$ be a Maass cusp form, with Laplace eigenvalue $\lambda = \frac{1}{4} + R^2$, of level $N$ and trivial character. Moreover, let $a_f(n)$ be the Hecke eigenvalues of $f$. We define the associated $L$-function to $f$ by
\begin{align*}
	L_f(s) = \sum_{n=1} \frac{a_f(n)}{n^s},
\end{align*}
where $\Real(s) > 1$. This can be analytically continued to the whole complex plane and satisfies the functional equation
\begin{align*}
	\Lambda_f(s) = N^{\frac{s}{2}}\Gamma_{\RR}(s + a + iR) \Gamma_{\RR}(s + a - iR) L_f(s)  = \omega (-1)^a \Lambda_f(1-s),
\end{align*}
where
\begin{itemize}
\item $\Gamma_{\RR}(s) = \pi^{-s/2} \Gamma(s/2),$
\item $\omega$ is the eigenvalue of the Fricke involution given by $f(z) = \omega f \left( -\frac{1}{Nz} \right)$,
\item $a=0$ if $f$ is even and $a=1$ if $f$ is odd. 
\end{itemize}

\begin{figure}[ht]
    \centering
    \def\svgwidth{\linewidth}
	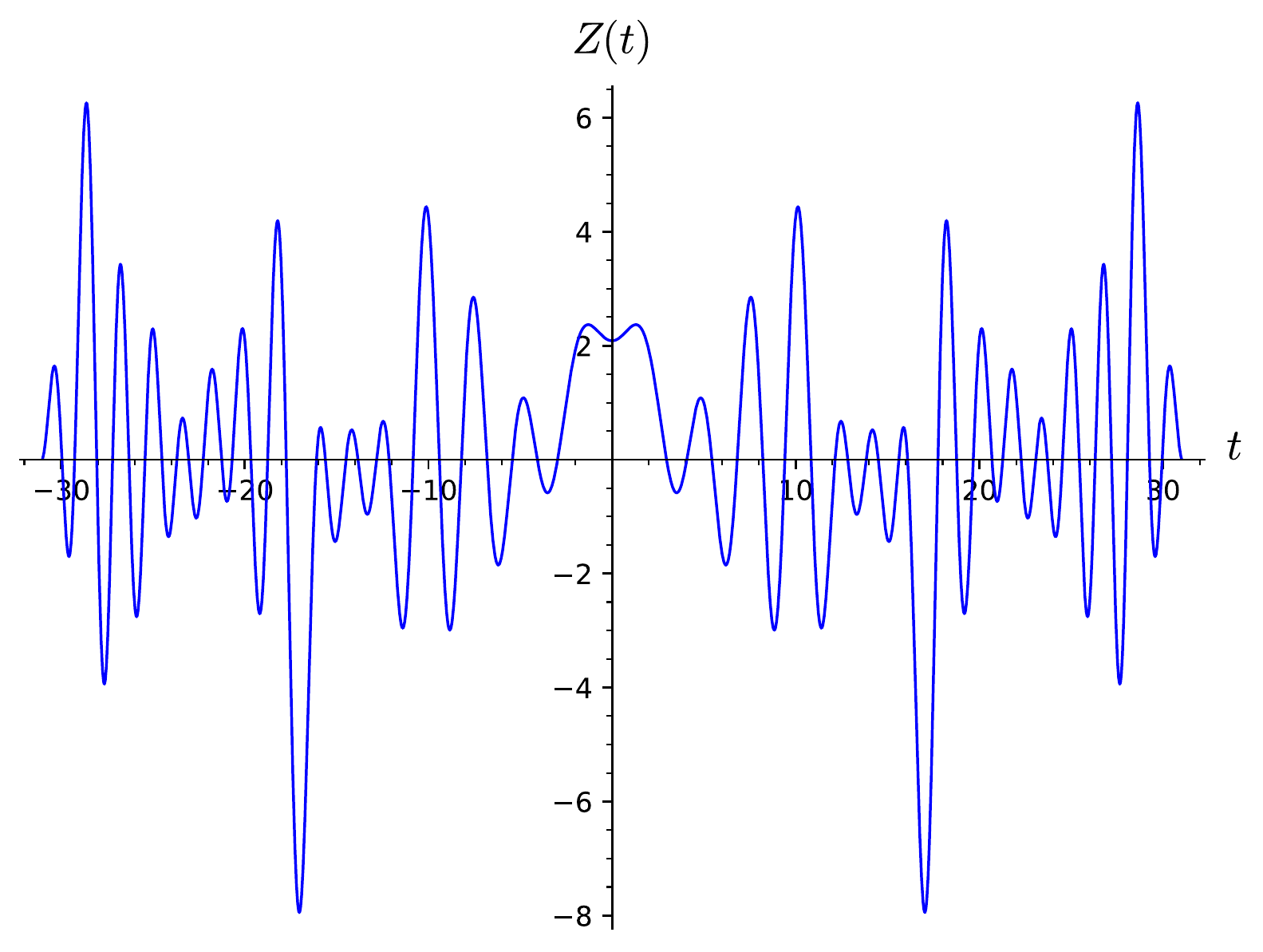
    \caption{Plot of the $Z$-function on the critical line associated to the first level $105$ Maass cusp form with Laplace eigenvalue $R = 0.4366582 \ldots$.}
    \label{fig:zfunc}
\end{figure}

It is conjectured, analogous to the Riemann zeta function, that $L$-functions associated to Maass cusp forms on $\Gamma_0(N)$ satisfy a Riemann hypothesis, that is all the zeros of $L_f(s)$ in the strip $\lbrace s \in \CC | 0 < \Real(s) < 1 \rbrace$ lie on the line $s = 1/2 + it, t \in \RR$. When computing zeros on the critical line of these $L$-functions, it is easier to work with the associated real-valued $Z$-function, defined by
\begin{align*}
    Z(t) = \bar{\varepsilon}^{1/2} \frac{\gamma(1/2+it)}{|\gamma(1/2+it)|} L_f(1/2+it),
\end{align*}
where $\gamma(s) = N^{\frac{s}{2}}\Gamma_{\RR}(s + a + iR) \Gamma_{\RR}(s + a - iR)$ and $\varepsilon = \omega (-1)^a$. Since $|Z(t)| = |L_f(1/2+it)|$, they share the same zeros on the critical line. An example of a $Z$-function is shown in Figure \ref{fig:zfunc}.

For the Maass forms we computed we used Rubenstein's library \textbf{lcalc} \cite{lcalc} to compute the $L$-function and calculate the zeros in the strip. We did this for all complete forms with $\varepsilon_i \leq 10^{-10}$ and found no zeros off the line, up to height $t=100$. To do this we computed the $a_f(n)$ with $(n,N) > 1$ up to $n \leq 2000$ using the method in Section \ref{subsec:hecke_eigenvals_coprime}. The method employed in \textbf{lcalc} to find zeros on the critical line is heuristic, however computing zeros on the critical line could be made rigorous with more work using the method in \cite{bookerthen}.

\subsection*{Availability of data}
The code generated as part of this work is available at \cite{Andrei2022}. The dataset of Maass forms is available at \cite{Andrei_Dataset_2022}.

\newpage
\bibliographystyle{abbrv}
\bibliography{tracepaper.bib}

\begin{thebibliography}{10}

\bibitem{AtkinLehner}
A.~O.~L. Atkin and J.~Lehner.
\newblock Hecke operators on {$\Gamma _{0}(m)$}.
\newblock {\em Math. Ann.}, 185:134--160, 1970.

\bibitem{modformspari}
K.~Belabas and H.~Cohen.
\newblock Modular forms in {P}ari/{GP}.
\newblock {\em Res. Math. Sci.}, 5(3):Paper No. 37, 19, 2018.

\bibitem{AndyMin}
A.~R. Booker and M.~Lee.
\newblock The {S}elberg trace formula as a {D}irichlet series.
\newblock {\em Forum Math.}, 29(3):519--542, 2017.

\bibitem{stfz}
A.~R. Booker and A.~Str\"{o}mbergsson.
\newblock Numerical computations with the trace formula and the {S}elberg
  eigenvalue conjecture.
\newblock {\em J. Reine Angew. Math.}, 607:113--161, 2007.

\bibitem{bookerstromvenka}
A.~R. Booker, A.~Str\"{o}mbergsson, and A.~Venkatesh.
\newblock Effective computation of {M}aass cusp forms.
\newblock {\em Int. Math. Res. Not.}, pages Art. ID 71281, 34, 2006.

\bibitem{bookerthen}
A.~R. Booker and H.~Then.
\newblock Rapid computation of {$L$}-functions attached to {M}aass forms.
\newblock {\em Int. J. Number Theory}, 14(5):1459--1485, 2018.

\bibitem{bump}
D.~Bump.
\newblock {\em Automorphic forms and representations}, volume~55 of {\em
  Cambridge Studies in Advanced Mathematics}.
\newblock Cambridge University Press, Cambridge, 1997.

\bibitem{Hejhal}
D.~A. Hejhal.
\newblock On eigenfunctions of the {L}aplacian for {H}ecke triangle groups.
\newblock In {\em Emerging applications of number theory ({M}inneapolis, {MN},
  1996)}, volume 109 of {\em IMA Vol. Math. Appl.}, pages 291--315. Springer,
  New York, 1999.

\bibitem{iwaniec}
H.~Iwaniec.
\newblock {\em Spectral methods of automorphic forms}, volume~53 of {\em
  Graduate Studies in Mathematics}.
\newblock American Mathematical Society, Providence, RI; Revista Matem\'{a}tica
  Iberoamericana, Madrid, second edition, 2002.

\bibitem{Jacobsonclassnumbers}
M.~J. Jacobson, Jr.
\newblock Experimental results on class groups of real quadratic fields
  (extended abstract).
\newblock In {\em Algorithmic number theory ({P}ortland, {OR}, 1998)}, volume
  1423 of {\em Lecture Notes in Comput. Sci.}, pages 463--474. Springer,
  Berlin, 1998.

\bibitem{arb}
F.~Johansson.
\newblock Arb: efficient arbitrary-precision midpoint-radius interval
  arithmetic.
\newblock {\em IEEE Transactions on Computers}, 66:1281--1292, 2017.

\bibitem{kimsarnak}
H.~Kim and S.~Peter.
\newblock Appendix 2: Refined estimates towards the ramanujan and selberg
  conjectures.
\newblock {\em Journal of the American Mathematical Society}, 16(1):175--181,
  2003.

\bibitem{LuoSarnakspacing}
W.~Luo and P.~Sarnak.
\newblock Number variance for arithmetic hyperbolic surfaces.
\newblock {\em Comm. Math. Phys.}, 161(2):419--432, 1994.

\bibitem{PARI}
{PARI~Group, The}, Univ. Bordeaux.
\newblock {\em {PARI/GP version \texttt{2.13.2}}}, 2020.
\newblock available from \url{http://pari.math.u-bordeaux.fr/}.

\bibitem{Risager}
M.~S. Risager.
\newblock Asymptotic densities of {M}aass newforms.
\newblock {\em J. Number Theory}, 109(1):96--119, 2004.

\bibitem{lcalc}
M.~O. Rubenstein.
\newblock {\em lcalc}.
\newblock {\tt http://code.google.com/p/l-calc/}.

\bibitem{sarnak}
P.~Sarnak.
\newblock Statistical properties of eigenvalues of the {H}ecke operators.
\newblock In {\em Analytic number theory and {D}iophantine problems
  ({S}tillwater, {OK}, 1984)}, volume~70 of {\em Progr. Math.}, pages 321--331.
  Birkh\"{a}user Boston, Boston, MA, 1987.

\bibitem{Andrei_Dataset_2022}
A.~Seymour-Howell.
\newblock {\em Dataset of Maass forms of squarefree level computed via the
  Trace Formula, Zenodo repository}, 2022.
\newblock \url{https://doi.org/10.5281/zenodo.7105772}.

\bibitem{Andrei2022}
A.~Seymour-Howell.
\newblock {\em Maass-Form-Trace-Formula-Code, GitHub repository}, 2022.
\newblock
  \url{https://www.github.com/aseymourhowell/Maass-Form-Trace-Formula-Code}.

\bibitem{Steil}
G.~Steil.
\newblock Eigenvalues of the {L}aplacian and of the {H}ecke operators for
  $\text{PSL}(2, \mathbb{Z})$.
\newblock {\em DESY}, (94--28), Hamburg 1994.

\bibitem{stromberg}
F.~Str\"omberg.
\newblock {\em Computational aspects of Maass Waveforms}.
\newblock PhD thesis, Uppsala University, 2005.

\bibitem{Andreaspre}
A.~Str\"ombergsson.
\newblock Explicit trace formula for {Hecke} operators.
\newblock {\em Preprint}, 2016.

\bibitem{sagemath}
{The Sage Developers}.
\newblock {\em {S}ageMath, the {S}age {M}athematics {S}oftware {S}ystem
  ({V}ersion 9.1.0)}, 2020.
\newblock {\tt https://www.sagemath.org}.

\end{thebibliography}

\end{document}